\documentclass[a4paper,leqno]{amsart}
\usepackage{amsfonts}
\theoremstyle{plain}
\begingroup

\newtheorem{theorem}{Theorem}[section]
\newtheorem{lemma}[theorem]{Lemma}
\newtheorem{proposition}[theorem]{Proposition}

\newtheorem{remark}[theorem]{Remark}
\newtheorem{definition}[theorem]{Definition}
\newtheorem{example}[theorem]{Example}
\endgroup

\theoremstyle{definition}
\theoremstyle{remark}

\mathchardef\emptyset="001F

\numberwithin{equation}{section}

\newcommand{\e}{\varepsilon}
\newcommand{\Om}{\Omega}
\newcommand{\Omk}{{\Om\setmeno K}}

\newcommand{\R}{{\mathbb R}}
\newcommand{\wto}{{\rightharpoonup}}
\newcommand{\setmeno}{\!\setminus\!}
\newcommand{\pinfty}{{+}\infty}

\newcommand{\huno}{{\mathcal H}^{1}}
\newcommand{\E}{{\mathcal E}}
\newcommand{\hb}{L^{1,2}}
\newcommand{\gdot}{{\dot{g}}}
\newcommand{\K}{{\mathcal K}(\overline\Om)}

\newcommand{\Km}{{\mathcal K}_m(\overline\Om)}
\newcommand{\Kmf}{{\mathcal K}_m^f(\overline\Om)}
\newcommand{\Kml}{{\mathcal K}_m^\lambda(\overline\Om)}

\title[Neumann problems in domains with cracks]
{Solutions of Neumann problems in domains with cracks\\
and applications to fracture mechanics}
\author[Gianni Dal Maso]{Gianni Dal Maso}
\thanks{Lecture notes of a course held in the 2001 CNA Summer School
``Multiscale Problems in Nonlinear Analysis'',
Carnegie Mellon University, Pittsburgh, May 31--June 9, 2001.}
\address{Gianni Dal Maso, SISSA, Via Beirut 2-4, 34014 Trieste,
Italy}
\email{dalmaso@sissa.it}
\begin{document}

\begin{abstract}
The first part of the course is devoted to the study of solutions to
the Laplace equation in $\Omega\setminus K$, where $\Omega$ is
a two-dimensional smooth domain and $K$ is a compact one-dimensional
subset of $\Omega$. The solutions are required to satisfy a homogeneous
Neumann boundary condition on $K$ and a nonhomogeneous Dirichlet
condition on (part of) $\partial\Omega$. The main result is the
continuous dependence of the solution on $K$, with respect to the
Hausdorff metric, provided that the number of connected components
of $K$ remains bounded. Classical examples show that the result is no
longer true without this hypothesis.

Using this stability result, the second part of the course develops
a rigorous mathematical formulation of a variational quasi-static
model of the slow growth of brittle fractures, recently introduced
by Francfort and Marigo. Starting from a discrete-time formulation,
a more satisfactory continuous-time formulation is
obtained, with full justification of the convergence arguments.
\end{abstract}
\maketitle

{\small


\bigskip
}

\bigskip

\begin{section}{NEUMANN PROBLEMS IN DOMAINS WITH CRACKS}

In these lectures $\Om$ is a fixed {\it bounded connected open\/}
subset of $\R^{2}$ with a {\it Lipschitz boundary\/} $\partial\Om$,
and $\K$ is the set of all compact subsets of $\overline\Om$.
Given $K\in\K$, we consider the differential equation $\Delta u=0$ on
${\Om\setmeno K}$ with a homogeneous Neumann boundary
condition on $\partial K$ and
on a part $\partial_N\Om$ of the boundary of $\Om$, and with a
non-homogeneous Dirichlet boundary condition on the rest of the boundary of
${\Om\setmeno K}$. For simplicity we assume that $\partial_N\Om$  is
a (possibly empty)
relatively open subset of $\partial\Om$, with a finite
number of connected components, and we set 
$\partial_D\Om:=\partial\Om\setmeno \overline{\partial_N\Om}$,
which turns out to be a relatively open subset of $\partial\Om$,
with a finite number of connected components. For the applications we
have in mind $K$ will be a one-dimensional set, and will be considered
as a crack in the domain $\Om$, but we do not need
this assumption in this lecture.

If $\partial K$ is not smooth, the
variational formulation of this problems requires the {\it Deny-Lions space\/}
$\hb(A)$, defined for every open set $A\subset\R^2$ as the space of
functions $u\in L^2_{loc}(A)$ such that the distributional gradient
$\nabla u$ belongs to $L^2(A;\R^2)$ (see \cite{DenLio}). The advantage of
this space is that the set $\{\nabla u: u\in \hb(A)\}$ is closed in
$L^2(A;\R^2)$ even if $\partial A$ is irregular
(for the proof  we refer, e.g., to \cite[Section 1.1.13]{Ma}). It is
well known that, if $A$ is bounded and has Lipschitz boundary, then
$L^2(A;\R^2)$ coincides with the usual Sobolev space $H^1(A)$ (see,
e.g., \cite[Corollary to Lemma 1.1.11]{Ma}).

Given  $K\in\K$ and $g\in \hb(\Omk)$,
we consider the following boundary value problem:
\begin{equation}\label{*}
\left\{\begin{array}{ll}
\Delta u=0 & \hbox{in }\Omk\,,\\
\noalign{\vskip5pt}
\displaystyle
\frac{\partial u}{\partial\nu}=0 & \hbox{on }
\partial(\Omk)\cap(K\cup \partial_N\Om)\,,\\
\noalign{\vskip5pt}
u=g & \hbox{on }\partial_D\Om\setmeno K\,.
\end{array}\right.
\end{equation}
By a solution of (\ref{*}) we mean a function $u$ which satisfies the
following conditions:
\begin{equation}\label{**}
\left\{\begin{array}{l}
u\in \hb(\Omk)\,,\ \, u=g\,
\hbox{ on }\,\partial_D\Om\setmeno K\,, \\
\noalign{\vskip5pt}
\displaystyle\int_{\Om\setminus K}\nabla u\,\nabla z\,dx=0
\quad\forall z\in
\hb(\Omk)\,,\ z=0\,\hbox{ on }\,
\partial_D\Om\setmeno K\,,
\end{array}\right.
\end{equation}
where the equalities on $\partial_D\Om$ are in the sense of traces.

It is clear that problem (\ref{**}) can be solved separately in
each connected component of $\Omk$. By
the Lax-Milgram lemma there exists
a unique solution in those components whose boundary meets
$\partial_D\Om\setmeno K$, while on the other components the solution is
given by an arbitrary constant. Thus the solution is not unique, if
there is a connected component whose boundary does not meet
$\partial_D\Om\setmeno K$. Note, however, that $\nabla u$ is
always unique.
Moreover, the map $g\mapsto\nabla u$ is linear from $\hb(\Omk)$
into $L^2(\Omk;\R^2)$ and satisfies the estimate
\begin{equation}\label{standard}
\int_{\Om\setminus K}|\nabla u|^2\,dx\le
\int_{\Om\setminus K}|\nabla g|^2\,dx\,.
\end{equation}

By standard arguments on the minimization of quadratic forms it is
easy to see that $u$ is a solution of problem (\ref{**}) if and only 
if $u$ is a
solution of the minimum problem
\begin{equation}\label{minpb}
\min_{v}\Big \{\int_{\Om\setminus K}|\nabla v|^2\,dx:
v\in\hb(\Om\setmeno K)\,,\ \, v=g\, \hbox{ on }\,
\partial_D\Om\setmeno K\Big\}\,.
\end{equation}
In these lectures, given a function $u\in \hb(\Omk)$ for some
$K\in\K$, we always
extend $\nabla u$ to $\Om$ by setting $\nabla u=0$ a.e.\ on $K$. Note that,
however, $\nabla u$ is the  distributional gradient  of $u$
only in $\Omk$, and, in
general, it does not coincide in $\Om$ with the gradient
of an extension of $u$.

To study the continuous dependence on $K$ of the solutions of problem
(\ref{**}) we consider the {\it Hausdorff distance\/} between two sets
$K_{1},\, K_{2}\in\K$, which is defined by
$$
d_{H}(K_{1},K_{2}):=
\max\big\{ \sup_{x\in K_1}{\rm dist}(x,K_2),
\sup_{y\in K_2}{\rm dist}(y,K_1)\big\}\,,
$$
with the conventions  ${\rm dist}(x,\emptyset)={\rm diam}(\Om)$ and
$\sup\emptyset=0$, so that $d_{H}(\emptyset, K)=0$ if $K=\emptyset$
and $d_{H}(\emptyset, K)={\rm diam}(\Om)$ if $K\neq\emptyset$.
We say that $(K_n)$ {\it converges to $K$ in the Hausdorff metric\/} if
${d_H(K_n,K)\to0}$. The following compactness theorem is
well-known (see, e.g., 
\cite[Blaschke's Selection Theorem]{Rog}).
\begin{theorem}\label{compactness}
Let $(K_n)$ be a sequence in $\K$.
Then there exists a subsequence which
converges in the Hausdorff metric to a set
$K\in \K$.
\end{theorem}

The following example shows that the convergence of
$(K_n)$ to $K$ in the Hausdorff metric does not imply the convergence
of the solutions of problems (\ref{**})
relative to $K_n$ to the solution relative to $K$, if we have no
bound on the number of connected components of $K_n$. In the next
lecture we will prove a convergence result under a uniform
bound on the number of connected components of~$K_n$.

\begin{example}\label{example}
Let $\Om:=(0,1)\times(-1,1)$, let $\partial_D\Om:=(0,1)\times\{-1,1\}$,
and let $g:=\pm 1$ on $(0,1)\times\{\pm 1\}$. For every $n\ge1$ let
$$
K_n:=\bigcup_{0=1}^{n-1}
\Big[\frac in, \frac in+\frac1{2n}\Big]\times\{0\}\,,
$$
and let $u_n$ be the solution of problem (\ref{*}) relative to $K_n$
and $g$.
Then $(K_n)$ converges in the Hausdorff metric to the set
$K:=[0,1]\times\{0\}$ and $(u_n)$ converges in $L^2(\Om)$
to the solution $u$ of the problem
\begin{equation}\label{*ex}
\left\{\begin{array}{ll}
\Delta u=0 & \hbox{in }\Om\,,\\
\noalign{\vskip5pt}
\displaystyle
\frac{\partial u}{\partial\nu}=0 & \hbox{on }
\partial_N \Om=\partial \Om\setmeno \overline{\partial_D \Om}\,,\\
\noalign{\vskip5pt}
u=g & \hbox{on }\partial_D\Om\,.\\
\end{array}\right.
\end{equation}
Since this solution is given explicitly by $u(x_1,x_2):=x_2$,
we see that it does not satisfy
the Neumann boundary condition ${\partial u}/{\partial\nu}=0$ on $K$.
\end{example}

\begin{proof}
Let $\Om^\pm:=\{x\in\Om:\pm x_2>0\}$.
By (\ref{standard}) we have $\int_{\Om^\pm} |\nabla u_n|^2dx\le
\int_{\Om} |\nabla g|^2dx$ for
every $n$, so that, passing to a subsequence, we may assume that
$(u_n)$ converges weakly in $H^1({\Om^+\cup\Om^-})$
  to a function $w\in H^1({\Om^+\cup\Om^-})$ such that $w=g$ on
  $\partial_D\Om$.
By symmetry the traces $u^\pm_n$ of $u_n$ from $\Om^\pm$ vanish on
${K\setmeno K_n}$. Let $\chi_n$ be the characteristic function of
${K\setmeno K_n}$. {}From the definition of $K_n$ it follows that
$(\chi_n)$ converges to $1/2$ weakly in $L^2(K,\huno)$, where
$\huno$ is the one-dimensional Hausdorff measure. Since the trace
operator is compact from $H^1(\Om^\pm)$ into $L^2(K,\huno)$,
the traces $u^\pm_n$ of $u_n$ from $\Om^\pm$ converge
to the corresponding traces $w^\pm$ of $w$ strongly in
$L^2(K,\huno)$. Therefore $(u^\pm_n\chi_n)$ converges to
$w^\pm/2$ weakly in
$L^1(K,\huno)$. As $u^\pm_n\chi_n=0$ on $K$ (recall that $u^\pm_n=0$ on
${K\setmeno K_n}$ and $\chi_n=0$ on $K_n$), we conclude that $w^\pm/2=0$
on $K$, therefore $w\in H^1(\Om)$.
By using the weak formulation (\ref{**}) we obtain
$$
\int_{\Om\setminus K_n} \nabla u_n\,\nabla z\,dx=0\qquad \forall
z\in H^1(\Om)\,,\ \,z=0 \,\hbox{ on }\, \partial_D\Om\,.
$$
Passing to the limit as $n\to\infty$ we obtain
$$
\int_{\Om} \nabla w\,\nabla z\,dx=0\qquad \forall
z\in H^1(\Om)\,,\ \,z=0 \,\hbox{ on }\, \partial_D\Om\,.
$$
This implies that $w$
coincides with the solution of problem (\ref{*ex}). By uniqueness,
the whole sequence $(u_n)$ converges to~$u$.
\end{proof}

\begin{remark}{\rm
The hypothesis $g:=\pm 1$ on $[0,1]\times\{\pm 1\}$ was intoduced
only to simplify the proof. Indeed
the same result holds when $g$ is an arbitrary function of $H^1(\Om)$.
To prove this fact, we can not use the equality $u_n=0$ on
${K\setmeno K_n}$,
which is not true in the general case; instead we introduce the functions
$\tilde u_n(x_1,x_2):=u_n(x_1,x_2)-u_n(x_1,-x_2)$ and
$\tilde w(x_1,x_2):=w(x_1,x_2)-w(x_1,-x_2)$, and observe that
$(\tilde u_n)$ converges to $\tilde w$ weakly in $H^1({\Om^+\cup\Om^-})$.
Since the traces of $\tilde u_n$ from $\Om^\pm$ vanish on
${K\setmeno K_n}$, arguing as before we obtain that
the traces of $\tilde w$ vanish on $K$. This implies that $w^+=w^-$
on $K$, and hence $w\in H^1(\Om)$. The conclusion follows now as
as in the previous proof.
}\end{remark}

In some cases the limit problem can contain a transmission condition,
as shown by the following example, for which we refer to \cite{Dam}
and~\cite{Mur}. Note that in this case the one-dimensional measure
of $K_n$ converges
to the one dimensional measure of~$K$.

\begin{example}\label{example2}
Let $\Om:=(0,1)\times(-1,1)$, let $\partial_D\Om:=(0,1)\times\{-1,1\}$,
and let $g\in H^1(\Om)$. For every $n\ge1$ let
$$
K_n:=\bigcup_{i=0}^{n-1}
\Big[\frac in, \frac{i+1}n-e^{-n}\Big]\times\{0\}\,,
$$
and let $u_n$ be the solution of problem (\ref{*}) relative to $K_n$.
Then $(K_n)$ converges in the Hausdorff metric to the set
$K:=[0,1]\times\{0\}$ and $(u_n)$ converges in $L^2(\Om)$
to the solution $u$ of the problem
\begin{equation}\label{*ex2}
\left\{\begin{array}{ll}
\Delta u=0 & \hbox{in }\Omk\,,\\
\noalign{\vskip5pt}
\displaystyle
\frac{\partial u}{\partial\nu}=0 & \hbox{on }
\partial \Om\setmeno (\overline{\partial_D\Om}\cup K) \,,\\
\noalign{\vskip5pt}
u=g & \hbox{on }\partial_D\Om\,,\\
\noalign{\vskip5pt}
\displaystyle
\frac{\partial u^\pm}{\partial\nu^\pm}=\pm \frac\pi2(u^- -u^+)&\hbox{on } K\,,
\end{array}\right.
\end{equation}
where $u^\pm$ is the restriction of $u$ to
$\Om^\pm:=\{x\in\Om:\pm x_2>0\}$ and $\nu^\pm$ is the outer unit
normal to $\Om^\pm$.
\end{example}

In the literature
of homogenization theory one can find other examples where
the convergence in the Hausdorff
metric of $K_n$ to  $K$
does not imply the convergence of the solutions of the
Neumann problems on $\Omk_n$ to the solution of the
Neumann problem on $\Omk$ (see, e.g.,
\cite{Khr},  \cite{Att-Pic}, and \cite{Cor}). These
papers show also that the bound on the number of connected components
of $K_n$, that we shall consider in the next lecture, would not
be enough in dimension larger than two.

The one-dimensional Hausdorff measure $\huno$ is not lower
semicontinuous with respect to the convergence in the Hausdorff metric.
For instance, in Example~\ref{example} we have $\huno(K_n)=1/2$ for
every $n$, while for the limit set we have $\huno(K)=1$.
However the lower semicontinuity holds if we have a uniform bound on
the number of connected components, as shown by the following theorem.

\begin{theorem}\label{Golab2}
Let $(K_n)$ be a sequence  in
$\K$ which converges to $K$ in the Hausdorff metric. Assume that each
set $K_n$ has at most $m$ connected components.
Then $K$ has at most $m$ connected components and
$$
\huno(K\cap U)\le \liminf_{n\to\infty} \,\huno(K_n\cap U)
$$
for every open set $U\subset \R^2$.
\end{theorem}

\begin{proof}
The case $m=1$ and $U=\R^2$ is the Go\l \c ab theorem,
for which we refer to \cite[Theorem~3.18]{Fal}. For an independent proof for an
arbitrary $U$ see \cite[Theorem 10.19]{MS}. The case $m>1$ is an easy
consequence of the case $m=1$ (see \cite[Corollary~3.3]{DM-Toa}).
\end{proof}
\end{section}

\begin{section}{CONVERGENCE OF SOLUTIONS}

Given an integer $m\ge 1$ and a constant $\lambda\ge 0$, let
$\Kml$ be the set of all compact
subsets $K$ of $\overline\Om$, with
$\huno(K)\le\lambda$, having at most $m$ connected components.

In this lecture we give the main ideas of the proof of the following
theorem. The complete proof can be found in \cite[Section~5]{DM-Toa}.

\begin{theorem}\label{convsol}
Let $m\ge1$ and $\lambda\ge0$, let $(K_n)$ be a sequence
in $\Kml$ which converges to $K$ in the Hausdorff metric, and let
$(g_n)$ be a sequence in $H^1(\Om)$ which converges to
$g$ strongly in $H^1(\Om)$.
Let $u_n$ be a solution
of the minimum problem
\begin{equation}\label{Pn}
\min_{v} \Big\{
\int_{\Om\setminus K_{n}}|\nabla v|^2\,dx:
v\in\hb(\Om\setmeno K_n)\,,\,\  v=g_n\, \hbox{ on }\,
\partial_D\Om\setmeno K_n \Big\}
\,,
\end{equation}
and let $u$ be  a solution of the minimum problem
\begin{equation}\label{P}
\min_{v} \Big\{\int_{\Om\setminus K}|\nabla v|^2\,dx:
v\in\hb(\Om\setmeno K)\,,\,\ v=g\, \hbox{ on }\,
\partial_D\Om\setmeno K \Big\}\,.
\end{equation}
Then $\nabla u_n\to \nabla u$ strongly in $L^2(\Om;\R^2)$.
\end{theorem}

This result is related to those
obtained by A. Chambolle and F. Doveri in \cite{Ch-D} and by D. Bucur
and N. Varchon in \cite{BucVar1}, \cite{BucVar2}, and \cite{BucVar},
which deal with the case of a pure Neumann
boundary condition. Since we impose a Dirichlet boundary condition
on $\partial_D\Om\setmeno K_n$ and a Neumann
boundary condition on the
rest of the boundary, our results can not be deduced easily from
these papers, so we give an independent proof, which uses the
duality argument which appears also in \cite{BucVar}.

To focus on the main ideas of the proof,
we consider only the case
$m=1$ and $\Om$ simply connected. Moreover we assume that
$K\cap\partial\Om=\emptyset$, to avoid minor difficulties arising at
the boundary. The technicalities needed to avoid
these simplifying hypotheses can be found in~\cite{DM-Toa}.

First of all, we want to construct the harmonic conjugate of $u_n$.
Let $R$ be the rotation on $\R^2$ defined by
$R(y_{1},y_{2}):=(-y_{2},y_{1})$.

\begin{definition}\label{constonK}
We say that a function $v\in H^1(\Om)$ is equal to a constant $c$
on a set $K\in\K$
if there exists a sequence $(v_n)$ in $C^1(\overline\Om)$
converging to $v$ strongly in $H^1(\Om)$
and such that each $v_n$ is equal $c$ in
a neighbourhood of~$K$.
\end{definition}

\begin{remark}\label{CoK}{\rm
It is possible to prove that, if
$$
\lim_{r\to0}\,\mathop{\rm ess\,sup}_{B_r(x)} |v-c|=0\qquad
\forall x\in K\,,
$$
then $v=c$ on $K$ in the sense of the previous definition (it is
enough to adapt the proof of \cite[Th\'eor\`eme IX.17]{Brez} or to
apply \cite[Theorem~4.5]{HKM}).
}\end{remark}

\begin{theorem}\label{conjug} Let $K$ be a connected compact set
contained in $\Om$ and
let $u$ be a solution of problem (\ref{**}). Then
there exists a function $v\in H^1(\Om)$ such that
$\nabla v=R\,\nabla u$ a.e.\ on $\Om$. Moreover
$v$ is constant on $K$ and on each connected component of
$\partial_N\Om$ (according to Definition~\ref{constonK}).
\end{theorem}

\begin{proof} If $\varphi\in C^\infty_c(\Om)$, we have
\begin{equation}\label{div}
\int_{\Om} \nabla u\,\nabla\varphi\,dx =
\int_{\Om\setminus K} \nabla u\,\nabla\varphi\,dx=0\,,
\end{equation}
where the first equality follows from our convention $\nabla u=0$
a.e.\ in $K$, while the second one follows from (\ref{**}).
Equality (\ref{div}) means that
${\rm div}(\nabla u)=0$ in ${\mathcal D}'({\Om})$, hence
${\rm rot}(R\,\nabla u)=0$ in ${\mathcal D}'({\Om})$. As
${\Om}$ is simply connected and has a Lipschitz boundary,
there exists $v\in H^1({\Om})$ such that $\nabla v=R\,\nabla u$
a.e.\ in ${\Om}$.

Since $\partial u/\partial \nu=0$ on $\partial_N\Om$, the tangential
derivative of $v$ (which is equal to the normal derivative of $u$)
vanishes on $\partial_N\Om$, and this implies that
$v$ is constant on each connected component of $\partial_N\Om$.

If $K$ has a non-empty interior and a smooth boundary, then $v$ is constant
a.e. on the interior of $K$, since $\nabla v=0$ a.e.\ on ${K}$.
Therefore $v$ is constant on $K$ according to Definition~\ref{constonK}.

The case of a general $K\in \K$ can be obtained by
approximating $K$ by a decreasing sequence of compact sets with
non-empty interior and a smooth boundary
(we refer to \cite[Theorem~4.2]{DM-Toa} for the details).
\end{proof}

\begin{proof}[Proof of Theorem~\ref{convsol}.]
Note that $u$ is a minimum point  of (\ref{P}) if and only if it satisfies
(\ref{**}); similarly, $u_n$ is a minimum point  of (\ref{Pn})
if and only if it satisfies (\ref{**}) with $K$ and $g$ replaced
by $K_n$ and $g_n$.

Taking $v:=g_n$ in the functional to be minimized,
we obtain that the sequence
$(\nabla u_n)$ is bounded in $L^2(\Om;\R^2)$. Therefore,
passing to a subsequence, $(\nabla u_n)$ converges weakly
in $L^2(\Om;\R^2)$ to a function $\psi$. As
$K\in\Kml$ by Theorem~\ref{Golab2}, and hence ${\rm meas}(K)=0$,
is easy to see that there exists a function $u^*\in \hb(\Omk)$,
with $u^*=g$ on
${\partial_D\Om}$,
such that $\psi=\nabla u^*$ a.e.\ on $\Om$ (we are assuming here
that $K\cap \partial_D\Om=\emptyset$, see
\cite[Lemma~4.1]{DM-Toa} for the details).

We will prove that
\begin{equation}\label{2uu*}
\nabla u^*=\nabla u\,\hbox{ a.e.\ in }\Om\setmeno K\,.
\end{equation}
As the limit does not depend on the subsequence, this implies that the whole
sequence $(\nabla u_n)$ converges to $\nabla u$ weakly in
$L^2(\Om;\R^2)$. Taking $u_n-g_n$ and $u-g$ as test functions
in the equations satisfied
by $u_n$ and $u$, we obtain
$$
\int_{\Om}|\nabla u_n|^2dx=\int_{\Om}\nabla u_n \nabla g_n\,dx\,,
\qquad
\int_{\Om}|\nabla u|^2dx=\int_{\Om}\nabla u \,\nabla g\,dx\,.
$$
As $\nabla u_n \wto \nabla u$ weakly in $L^2(\Om;\R^2)$ and $\nabla
g_n\to\nabla g$ strongly  in $L^2(\Om;\R^2)$, from the
previous equalities we obtain that $\|\nabla u_n\|_{L^2(\Om;\R^2)}$
converges to $\|\nabla u\|_{L^2(\Om;\R^2)}$, which implies the strong
convergence of the gradients in $L^2(\Om;\R^2)$.

By the uniqueness of the gradients of the solutions,
to prove (\ref{2uu*}) it is enough to show that $u^*$ is a solution of
(\ref{**}). This will be done by considering, for each $u_n$, its
harmonic conjugate $v_n$ given by Theorem~\ref{conjug}.
By adding a suitable constant,
we may assume that $\int_{\Om} v_n\,dx=0$ for every $n$.
Since $\nabla v_n=R\,\nabla u_n$ a.e.\ on ${\Om}$,
we deduce that $(\nabla v_n)$ converges to $R\,\nabla u^*$ weakly in
$L^2(\Om;\R^2)$, and by the Poincar\'e inequality
$(v_n)$  converges weakly in $H^1(\Om)$ to a
function $v$ which satisfies
$\nabla v=R\,\nabla u^*$ a.e.\ on~$\Om$.

Let us prove that $v$ is constant on $K$ according
to Definition~\ref{constonK}. This is trivial if $K$
reduces to a point. If $K$ has more than one point,
then $\liminf_{n}{\rm diam}(K_n)>0$;
since the sets $K_n$ are connected, we obtain also
$\liminf_{n}{\rm cap}(K_n,\Om)>0$, where the capacity
${\rm cap}(K_n,\Om)$ of $K_n$ with respect to $\Om$ is defined by
$$
{\rm cap}(K_n,\Om):=\inf_z \Big\{\int_\Om |\nabla z|^2 dx:
z\in C^1_c(\Om)\,,\ \,z\ge 1 \, \hbox{ on }\, K_n\Big\}\,.
$$
As $v_n=c_n$ on $K_n$ for
suitable constants $c_n$ (see Theorem~\ref{conjug}),
using the Poincar\'e inequality (see, e.g., \cite[Corollary 4.5.3]{Zie})
it follows that $(v_n-c_n)$ is bounded in $H^1(\Om)$, hence
the sequence $(c_n)$ is bounded, and therefore, passing to a
subsequence, we may assume that $(c_n)$ converges to a suitable
constant $c$.

Let us fix a constant $R>0$ with $R<{\rm diam}(K)/2$ and
$R<{\rm dist}(K,\partial\Om)$.
Since $\Delta v_n=0$ on ${\Om\setmeno K_n}$,
by Maz'ya's estimate (see \cite[Theorem~1]{Maz67}) there exist two
constants $M>0$ and $\beta>0$, independent of $n$, $x_n$, and $r$, such that
for every $x_n\in K_n$ and every $r\in(0,R)$
\begin{equation}\label{Mazya}
\mathop{\rm ess\,sup}_{B_r(x_n)}
|v_n-c_n|
\le M \exp\Big(-\beta\int_r^R
\gamma_n(x_n,\rho)\,\frac{d\rho}{\rho} \,\Big)\,,
\end{equation}
where
$\gamma_n(x_n,\rho):={\rm cap}(K_n\cap B_\rho(x_n),B_{2\rho}(x_n))$.
For $n$ large we have ${\rm diam}(K_n)>2R$, so that for every $x_n\in
K_n$ and every $\rho\in(0,R)$ we have
$K_n\cap \partial B_\rho(x_n)\neq\emptyset$. As $K_n$ is connected,
there exists a constant $\alpha>0$, independent of $n$, $x_n$, and $\rho$,
such
that
$\gamma_n(x_n,\rho)= {\rm cap}(K_n\cap B_\rho(x_n),B_{2\rho}(x_n))\ge
\alpha$. Therefore (\ref{Mazya}) yields
$$
\mathop{\rm ess\,sup}_{B_r(x_n)}
|v_n-c|
\le M\,\Big(\frac{r}{R} \Big)^{\alpha \beta} + |c_n-c|
$$
for every $x_n\in K_n$ and every $r\in(0,R)$.

Let us fix $x\in K$ and a sequence $x_n\in K_n$
converging to $x$. For every $\rho\in (0,r)$ we have $B_\rho(x)\subset
B_r(x_n)$ for $n$ large enough, hence
$$
\mathop{\rm ess\,sup}_{B_\rho (x)}
|v_n-c|
\le M\,\Big(\frac{r}{R} \Big)^{\alpha \beta} + |c_n-c|
$$
for $n$ large enough. Passing to the limit first as $n\to\infty$ and then
as $\rho\to r$ we get
\begin{equation}\label{Mazya2}
\mathop{\rm ess\,sup}_{B_r (x)}
|v-c|
\le M\,\Big(\frac{r}{R} \Big)^{\alpha \beta}\,.
\end{equation}
As $r\to0$ we obtain that $v$ is equal to $c$
on~$K$ (see Remark~\ref{CoK}).

On the other hand, every $v_n$ is constant on each
connected component of $\partial_N\Om$ (see Theorem~\ref{conjug}).
Since
$v_n\wto v$ weakly in $H^1(\Om)$, we conclude
that $v$ is constant on each connected component of
$\partial_N\Om$.

Therefore there exists a sequence $(w_n)$ in $C^1(\overline\Om)$
converging to $v$ strongly in $H^1(\Om)$,
such that each $w_n$ is constant in a neighbourhood of $K$ and
in a neighbourhood of each connected component of $\partial_N\Om$.

Let $z\in \hb(\Omk)$
with $z=0$ on $\partial_D\Om$., and let
$\varphi_n\in C^1(\overline\Om)$ with $\varphi_n=1$ on
${\rm supp}(\nabla w_n)$ and $\varphi_n=0$ on a neighbourhood
of $K\cup\partial_N\Om$.
As ${\rm div}(R\,\nabla w_n)=0$ in $\Om$ and
$z\,\varphi_n\in H^1_0(\Om)$, we
have
$$
\int_{\Om\setminus K}R\,\nabla w_n\,\nabla z\,dx=
\int_{\Om\setminus K}R\,\nabla w_n\,\nabla (z\,\varphi_n)\,dx
  =0\,.
$$
Since $R\,\nabla u^*=\nabla v$ a.e.\ on $\Om$,
passing to the limit as $n\to\infty$ we obtain
$$
\int_{\Om\setminus K}\nabla u^*\,\nabla z\,dx=
-\int_{\Om\setminus K}R\,\nabla v\,\nabla z\,dx=0\,,
$$
which shows that $u^*$ is a solution of (\ref{**}).
\end{proof}
\end{section}

\begin{section}{A QUASI-STATIC MODEL FOR BRITTLE FRACTURES}

Since the pioneering work of A.~Griffith \cite{Gri}, the growth of a
brittle fracture is considered to be the result of the competition
between the energy spent to increase the crack and the corresponding
release of bulk energy. This idea is the basis of the celebrated
Griffith's criterion for crack growth (see, e.g., \cite{SL}), and is used to
study the crack propagation along a preassigned path.
The actual path followed by
the crack is often determined by using different criteria
(see, e.g., \cite{ES}, \cite{SL}, \cite{SM}).

Recently G.A.~Francfort and J.-J.~Marigo \cite{FraMar3} proposed
a variational model for the quasi-static growth of brittle
fractures, based on Griffith's theory,
where the interplay between bulk and surface energy
determines also the crack path.

The purpose of this and of the next lecture
is to give a precise mathematical
formulation of a variant of this model in the {\it two-dimensional
case\/}, and to prove an existence result for the {\it quasi-static 
evolution of a
fracture\/} by using the {\it time discretization method\/} proposed
in~\cite{FraMar3}.

To simplify the mathematical description of the model, we consider
only {\it linearly elastic homogeneous isotropic
materials\/}, with Lam\'e coefficients
$\lambda$ and $\mu$. We restrict our analysis to
the case of an {\it anti-plane shear\/}, where
the reference configuration is an infinite cylinder ${\Om{\times}\R}$, with
$\Omega\subset \R^2$, and the displacement has the special form
$(0,0,u(x_1,x_2))$ for every $(x_1,x_2,y)\in {\Om{\times}\R}$.
We assume also that the cracks  have the form
$K{\times}\R$, where $K$ is a compact set in $\overline\Om$.
In this case the notions of bulk energy and surface energy
refer to a finite portion of the cylinder determined by
two cross sections separated by a unit distance.
The {\it bulk energy\/} is given by
\begin{equation}\label{bulk}
\frac{\mu}2\int_{\Om\setminus K} |\nabla u|^2dx\,,
\end{equation}
while the {\it surface energy\/} is given by
\begin{equation}\label{surface}
k\,\huno(K)\,,
\end{equation}
where $k$ is a constant which depends on the toughness of the material, and
$\huno$ is the {\it one-dimensional Hausdorff measure\/}, which
coincides with the ordinary length in case $K$ is a rectifiable arc. For
simplicity we take $\mu=2$ and $k=1$ in (\ref{bulk})
and (\ref{surface}).

As in the previous lectures $\Om$ is a fixed {\it bounded connected open\/}
subset of $\R^{2}$ with {\it Lipschitz boundary\/}. Following
\cite{FraMar3}, we fix an open subset $\partial_D\Om$ of $\partial\Om$, on
which we want to  prescribe a {\it Dirichlet boundary condition\/} for
the displacement $u$.
As in the previous lectures we assume that $\partial_D\Om$
has a {\it finite number of connected
components\/}.

Given a function $g$ on $\partial_D\Om$,
we consider the boundary condition
$u=g$ on $\partial_D\Omk$. We can not prescribe a Dirichlet
boundary condition on $\partial_D\Om\cap K$, because
the boundary displacement is not
transmitted through the crack, if the crack touches the boundary. Assuming
that {\it the fracture is traction free\/}
(and, in particular, without friction),
the displacement $u$ in $\Omk$ is obtained by
{\it minimizing (\ref{bulk}) under the boundary condition
$u=g$ on $\partial_D\Omk$\/}.
The {\it total energy\/} relative to the boundary displacement
$g$ and to the crack determined by $K$ is therefore
\begin{equation}\label{e}
\E(g,K)=\min_v \Big\{\int_{\Om\setminus K}|\nabla
v|^2dx+\huno(K) : v\in\hb(\Om),\ v=g \hbox{ on }\partial_D\Om\setmeno 
K \Big\}\,.
\end{equation}
The existence of a minimizer has been proved in the first lecture.

In the theory developed in \cite{FraMar3} a crack with finite
surface energy is any compact subset $K$ of $\overline\Om$ with
$\huno(K)<\pinfty$. For technical reasons, due to the hypotheses of
Theorem~\ref{convsol}, we propose a variant of this model, where we 
prescribe an a
priori bound on the number of connected components of the
cracks. Without this restriction, some
convergence arguments used in the proof of our existence result
are not justified
by the present development of
the mathematical theories related to this subject.
Given an integer $m\ge 1$, let $\Kmf$ be the set of all compact
subsets $K$ of $\overline\Om$, with $\huno(K)<\pinfty$,
having at most $m$ connected components.

We begin by describing a discrete-time model of {\it quasi-static irreversible
evolution of a fracture\/} under the action of a {\it time dependent
boundary displacement\/} $g(t)$, $0\le t\le 1$.
As usual, we assume that $g(t)$ can be extended to a function, still
denoted by $g(t)$, which belongs to the Sobolev space $H^1(\Omega)$.
For simplicity, we assume also that $g(0)=0$.

Given a time step $\delta>0$, for every integer $i\ge 0$ we set
$t_i^\delta:=i\delta$ and $g_i^\delta:=g(t_i^\delta)$.
The fracture $K_i^\delta$ at time $t_i^\delta$ is defined inductively
in the following way: for $i=0$ we set $K_0^\delta:=K_0$, while for
$i\ge1$ $K_i^\delta$ is any minimizer of the problem
\begin{equation}\label{pidelta}
\min_K\big\{ \E(g_i^\delta,K) : K \in \Kmf,\
K \supset K_{i-1}^\delta \big\}\,.
\end{equation}

\begin{lemma}\label{existspidelta}
There exists a solution of the minimum problem (\ref{pidelta}).
\end{lemma}
\begin{proof}
By definition $K_0^\delta:=K_0\in\Kmf$. Assume by induction that
$K_{i-1}^\delta\in \Kmf$ and let $\lambda$ be a constant such that
$\lambda>\E(g_i^\delta,K_{i-1}^\delta)$.
Consider a minimizing sequence $(K_n)$ of problem~(\ref{pidelta}).
We may assume that $K_n\in\Kml$ for every $n$.
By the Compactness Theorem~\ref{compactness}, passing to a subsequence,
we may assume that  $(K_n)$ converges in the Hausdorff metric  to
some  compact set $K$ containing $K_{i-1}^\delta$.
For every $n$ let $u_n$ be a
solution of the minimum  problem (\ref{e}) which defines
$\E(g_i^\delta,K_n)$.
By Theorem~\ref{convsol} $(\nabla u_n)$ converges
strongly in $L^2(\Om;\R^2)$ to $\nabla u$, where $u$ is a  solution of
the minimum problem (\ref{e}) which defines $\E(g_i^\delta,K)$.
By Theorem~\ref{Golab2} we have $K\in\Km$ and
$\huno(K)\leq\liminf_{n}\huno(K_n)\le\lambda$, hence $K\in\Kml$. As
$\|\nabla u\|=\lim_{n}\|\nabla u_n\|$,
we conclude that $\E(g_i^\delta,K)\leq\liminf_{n}\E(g_i^\delta,K_n)$.  Since
$(K_n)$ is a minimizing sequence, this proves that $K$  is a
solution of the minimum problem (\ref{pidelta}).
\end{proof}

Let $u_i^\delta$ be a solution of the minimum problem (\ref{e})
which defines $\E(g_i^\delta,K_i^\delta)$. Then the pair
$(u_i^\delta,K_i^\delta)$ minimizes the sum of the bulk and surface
energy among all $K \in \Kmf$ with
$K \supset K_{i-1}^\delta$ and among all
$u\in\hb({\Om\setmeno K})$
with $u=g$ on ${\partial_D\Om\setmeno K}$.

In order to pass to a continuous-time model, given $T>0$
we define the step
functions $K_\delta\colon[0,T]\to\Kmf$ and $u_\delta\colon
[0,T]\to L^2_{\rm loc}(\Om)$ by setting
\begin{equation}\label{kdt}
K_\delta(t):=K_{i}^\delta\qquad \hbox{and}\qquad
u_\delta(t):=u_{i}^\delta\qquad \hbox{for }
t_i^\delta\le t<t_{i+1}^\delta\,.
\end{equation}
Our purpose is to pass to the limit as $\delta\to 0$. To this aim we
use the following result, which is the analogue of the Helly theorem for
compact valued increasing functions.

\begin{theorem}\label{Helly}
Let $(K_n)$ be a sequence of increasing functions from $[0,T]$ into
$\K$, i.e., $K_n(s)\subset K_n(t)$ for $0\le s<t\le T$.
Then there exist a subsequence, still denoted by $(K_n)$, and an
increasing function $K\colon[0,T]\to\K$, such that $K_n(t)\to K(t)$ in
the Hausdorff metric for every $t\in[0,T]$.
\end{theorem}

To prove Theorem~\ref{Helly} we use the following result, which extends
another well known property of real valued monotone functions.

\begin{lemma}\label{diff3.5}
Let $K_1, \,K_2\colon [0,T]\to\K$ be two
increasing functions such that
\begin{equation}\label{K12}
K_1(s) \subset K_2(t)\qquad\hbox{and} \qquad
K_2(s) \subset K_1(t)
\end{equation}
for every $s$, $t\in[0,T]$ with $s<t$.
Let $\Theta$ be the set of points $t\in[0,T]$ such that  $K_1(t)=K_2(t)$.
Then ${[0,T]\setmeno\Theta}$ is at most countable.
\end{lemma}

\begin{proof}
For $i=1,2$, consider the functions
$f_i\colon \overline \Om\times[0,T]\to\R$ defined by
$f_i(x,t):={\rm dist}(x,K_i(t))$, with the convention
that ${\rm dist}(x,\emptyset)={\rm diam}(\Om)$. Then the
functions $f_i(\cdot,t)$ are
Lipschitz continuous with constant $1$ for every $t\in[0,T]$, and
the functions $f_i(x,\cdot)$ are non-increasing for every $x\in\overline\Om$.

Let $D$ be a countable dense subset of $\overline\Om$.  For every
$x\in D$ there exists a countable set $N_x\subset[0,T]$ such that
$f_i(x,\cdot)$ are continuous at every point of $[0,T]\setmeno N_x$.
By (\ref{K12}) we have $f_1(x,s)\ge f_2(x,t)$ and $f_2(x,s)\ge f_1(x,t)$ for
every $x\in\overline\Om$ and every $s$, $t\in[0,T]$ with $s<t$. This
implies that $f_1(x,t)=f_2(x,t)$ for every $x\in D$ and every
$t\in[0,T]\setmeno N_x$.
Let $N$ be the countable set defined by
$N:=\bigcup_{x\in D}N_x$, and let $t\in[0,T]\setmeno N$.
Then $f_1(x,t)=f_2(x,t)$ for every $x\in D$, and,
by continuity, for every
$x\in\overline\Om$, which yields $K_1(t)=K_2(t)$.
This proves that ${[0,T]\setmeno N}\subset\Theta$,
hence ${[0,T]\setmeno \Theta}\subset N$.
\end{proof}

\begin{proof}[Proof of Theorem~\ref{Helly}.]
Let $D$ be a countable dense subset of $(0,T)$. Using the Compactness
Theorem~\ref{compactness} and a diagonal
argument, we find a subsequence,  still denoted by $(K_n)$, and an
increasing function $K\colon D\to\K$, such that $K_n(t)\to K(t)$ in
the Hausdorff metric for every $t\in D$. Let
$K^-\colon (0,T]\to\K$ and $K^+\colon [0,T)\to\K$ be the increasing
functions defined by
\begin{eqnarray*}
& K^-(t):=\displaystyle{\rm cl}
\Big(\bigcup_{s<t,\; s\in D}K(s)\,\Big)
\qquad\hbox{for }0<t\leq T\,,\\
& K^+(t):=\displaystyle\bigcap_{s>t,\;s\in D }K(s)\qquad\hbox{for }0
\leq t<T\,,
\end{eqnarray*}
where ${\rm cl}$ denotes the closure. Let $\Theta$ be the set of points
$t\in[0,T]$ such that  $K^-(t)=K^+(t)$. As $K^-$ and $K^+$ satisfy
(\ref{K12}), by Lemma~\ref{diff3.5} the set
${[0,T]\setmeno\Theta}$ is at most countable.

Since $K^-(t)\subset K(t)\subset K^+(t)$ for every $t\in D$, we have
$K(t)=K^-(t)=K^+(t)$ for every $t\in {\Theta\cap D}$.
For every $t\in {\Theta\setmeno D}$ we define $K(t):=K^-(t)=K^+(t)$.
To prove that $K_n(t)\to K(t)$ for a given $t\in {\Theta\setmeno D}$,
by the Compactness Theorem~\ref{compactness} we may assume that
$K_n(t)$ converges in the Hausdorff metric to a set $K^*$. For
every $s_1,\, s_2\in D$, with $s_1<t<s_2$, by monotonicity we have
$K(s_1)\subset K^*\subset K(s_2)$. As $K^*$ is closed, this implies
$K^-(t)\subset K^*\subset K^+(t)$, therefore
$K_n(t)\to K(t)$ by the definitions of $\Theta$ and $K(t)$.

Since $[0,T]\setmeno(\Theta\cup D)$ is at most countable, by a
diagonal argument we find a further subsequence, still denoted by
$(K_n)$, and a function
$K\colon{[0,T]\setmeno (\Theta\cup D)}\to \K$, such that $K_n(t)\to
K(t)$  in the Hausdorff metric for every
$t\in {[0,T]\setmeno(\Theta\cup D)}$.

Therefore $K_n(t)\to K(t)$ for every $t\in[0,T]$, and this implies that
$K(\cdot)$ is increasing on $[0,T]$.
\end{proof}

The following result on the continuity of compact valued increasing maps
will be used in the next lecture. Its proof is similar to the proof of
Theorem~\ref{Helly}.

\begin{proposition}\label{diff3}
Let $K\colon [0,T]\to\K$ be an increasing function, and let
$K^-\colon (0,T]\to\K$ and $K^+\colon [0,T)\to\K$ be the functions
defined by
\begin{eqnarray}
&\textstyle K^-(t):={\rm cl}\big(\bigcup_{s<t}K(s)\big)
\qquad\hbox{for }0<t\leq T\,,
\label{k_*t}\\
&\textstyle K^+(t):=\bigcap_{s>t}K(s)\qquad\hbox{for }0
\leq t<T\,,\label{k^*t}
\end{eqnarray}
where ${\rm cl}$ denotes the closure. Then
\begin{equation}\label{7.9}
K^-(t)\subset K(t)\subset K^+(t)
\qquad\hbox{for }0< t<T\,.
\end{equation}
Let $\Theta$ be the set of points $t\in(0,T)$ such that  $K^-(t)=K^+(t)$.
Then ${[0,T]\setmeno\Theta}$ is at most countable, and $K(t_n)\to
K(t)$ in the Hausdorf metric for every $t\in\Theta$ and every sequence
$(t_n)$ in $[0,T]$ converging to~$t$.
\end{proposition}

\begin{proof}
It is clear that $K^+(\cdot)$ and $K^-(\cdot)$ are increasing and 
satisfy (\ref{K12}).
Therefore ${[0,T]\setmeno\Theta}$ is at most countable by
Lemma~\ref{diff3.5}.

Let us fix $t\in\Theta$ and a sequence $(t_n)$ in $[0,T]$ converging
to $t$. By the Compactness Theorem~\ref{compactness} we may assume
that $K(t_n)$ converges in the Hausdorff metric to a set $K^*$. For
every $s_1, \,s_2\in[0,T]$, with $s_1<t<s_2$, we have $K(s_1)\subset
K(t_n)\subset K(s_2)$ for $n$ large enough, hence $K(s_1)\subset
K^*\subset K(s_2)$. As $K^*$ is closed this implies $K^-(t)\subset
K^*\subset K^+(t)$, therefore $K^*=K(t)$ by (\ref{7.9}) and by the
definition of~$\Theta$.
\end{proof}

According to Theorem~\ref{Helly}, there exist a sequence
$(\delta_k)$ converging to $0$ and an increasing function
$K\colon[0,T]\to\K$ such that, for every $t\in[0,T]$, $K_\delta (t)\to K(t)$ in
the Hausdorff metric as $\delta$ tends to $0$ along this sequence.
In the next lecture we will prove the main properties of the function
$K\colon[0,T]\to\K$ obtained in this way, which represents the
continuous-time evolution of the fracture. To simplify
the notation, when we write $\delta\to0$ we always mean that
$\delta$ tends to $0$ along the sequence $(\delta_k)$ considered above.
\end{section}

\begin{section}{PROPERTIES OF THE CONTINUOUS-TIME MODEL}\label{irrev}

In this final lecture we prove some properties of the function
$K\colon[0,T]\to\K$ defined as the limit, as $\delta\to 0$, of the
functions $K_\delta\colon[0,T]\to\K$ given by the discrete-time model.
This function represents the quasi-static evolution of the crack in
our continuous-time model.

To prove these results, we assume that the function $t\mapsto g(t)$,
which gives the imposed boundary displacement on $\partial_D\Om$, is {\it
absolutely continuous\/} from $[0,T]$ into $H^1(\Omega)$. Its time
derivative is a Bochner
integrable function from $[0,T]$ into $H^1(\Omega)$, which will be
denoted by $\dot g(t)$. For the main properties of absolutely continuous
functions with values in a Hilbert space we refer, e.g.,
to \cite[Appendix]{Bre}.

We begin with a crucial estimate for the solutions of the
discrete-time problems. Here and in the rest of the lecture
$(\cdot|\cdot)$ and $\|\cdot\|$
denote the scalar product  and the norm in $L^2(\Om;\R^2)$.

\begin{lemma}\label{discr}
There exists a positive function $\omega(\delta)$, converging to zero
as $\delta\to0$, such that
\begin{equation}\label{2discr}
\|\nabla u_j^\delta\|^2+\huno(K_j^\delta)\leq
\|\nabla u_i^\delta\|^2+\huno(K_i^\delta)+
2\int_{t_i^{\delta}}^{t_j^{\delta}}(\nabla u_\delta(t)|\nabla
\gdot(t))\,dt+\omega(\delta)
\end{equation}
for $0\leq i<j$ with $t_j^{\delta}\le T$.
\end{lemma}
\begin{proof}
Let us fix an integer $r$ with $i\leq r<j$. {}From the absolute
continuity of $g$ we have
$$
g_{r+1}^\delta-g_r^\delta=\int_{t_r^{\delta}}^{t_{r+1}^{\delta}}\gdot(t)\,dt\,,
$$
where the integral is a Bochner integral for functions with values in
$H^1(\Om)$. This implies that
\begin{equation}\label{nabla}
\nabla g_{r+1}^\delta-\nabla g_r^\delta=\int_{t_r^{\delta}}^{t_{r+1}^{\delta}}
\nabla \gdot(t)\,dt\,,
\end{equation}
where the integral is a Bochner integral for functions with values in
$L^2(\Om;\R^2)$.

As $u_r^{\delta}+g_{r+1}^\delta-g_r^\delta\in \hb(\Omk_{r}^\delta)$
and $u_r^{\delta}+g_{r+1}^\delta-g_r^\delta=g_{r+1}^\delta$ on
$\partial_D\Omk_r^\delta$, we have
\begin{equation}\label{a}
\E(g_{r+1}^\delta,K_r^\delta)\le
\|\nabla u_r^\delta+\nabla g_{r+1}^\delta-\nabla g_r^\delta\|^2+
\huno(K_r^\delta)\,.
\end{equation}
By the minimality of $u_{r+1}^\delta$ and by (\ref{pidelta})
we have
\begin{equation}\label{b}
\|\nabla u_{r+1}^\delta\|^2+\huno(K_{r+1}^\delta)=
\E(g_{r+1}^\delta,K_{r+1}^\delta)\le
\E(g_{r+1}^\delta,K_r^\delta)\,.
\end{equation}
{}From (\ref{nabla}), (\ref{a}), and (\ref{b})  we obtain
\begin{eqnarray*}
& \|\nabla u_{r+1}^\delta\|^2+\huno(K_{r+1}^\delta)\leq
\|\nabla u_r^\delta+\nabla g_{r+1}^\delta-\nabla g_r^\delta\|^2+
\huno(K_r^\delta)\leq\\
&\leq\|\nabla u_{r}^\delta\|^2+\huno(K_r^\delta)+2
\displaystyle\int_{t_r^{\delta}}^{t_{r+1}^{\delta}}(\nabla u_{r}^\delta|\nabla
\gdot(t))\,dt+\Big(\displaystyle\int_{t_r^{\delta}}^{t_{r+1}^{\delta}}\|\nabla
\gdot(t)\|\,dt\Big)^2\leq\\
&\leq\|\nabla u_{r}^\delta\|^2+\huno(K_r^\delta)+2
\displaystyle\int_{t_r^{\delta}}^{t_{r+1}^{\delta}}
(\nabla u_\delta(t)|\nabla\gdot(t))\,dt+
\sigma(\delta)\displaystyle\int_{t_r^{\delta}}^{t_{r+1}^{\delta}}\|\nabla
\gdot(t)\|\,dt\,,
\end{eqnarray*}
where
$$
\sigma(\delta):=\max_{0\le r,\,t_r^{\delta}<T} 
\int_{t_r^{\delta}}^{t_{r+1}^{\delta}}\|\nabla
\gdot(t)\|\,dt \ \longrightarrow \ 0
$$
by the absolute continuity of the integral.
Iterating now this inequality for $i\leq r<j$ we get
(\ref{2discr}) with $\omega(\delta):=\sigma(\delta)\int_0^1\|\nabla
\gdot(t)\|\,dt$.
\end{proof}

\begin{lemma}\label{estim}
There exists a constant $\lambda$, depending only on $g$ and $K_0$,
such that
\begin{equation}\label{stima}
\|\nabla u_i^\delta\|\leq \lambda\quad\hbox{and}\quad
\huno(K_i^{\delta})\leq \lambda
\end{equation}
for every $\delta>0$ and for every  $i\ge 0$ with $t_i^\delta\le T$.
\end{lemma}
\begin{proof}
As $v:=g_i^\delta$ is admissible for the problem
(\ref{e}) which defines $\E(g_i^\delta, K_i^\delta)$,
by the minimality of $u_i^\delta$ we have
$\|\nabla u_i^\delta\|\leq\|\nabla g_i^\delta\|$, hence
$\|\nabla u_\delta(t)\|\leq\|\nabla g_\delta(t)\|$ for every
$t\in[0,T]$.
As $t\mapsto g(t)$ is absolutely continuous with values in $H^1(\Om)$,
the function $t\mapsto \|\nabla\gdot(t)\|$ is integrable on $[0,T]$
and
there exists a constant $C>0$ such that
$\|\nabla g(t)\|\leq C$ for every $t\in[0,T]$.  This implies the
former
inequality in (\ref{stima}). The latter inequality follows
now from Lemma~\ref{discr} and from the inequality
$\|\nabla u_0^\delta\|^2 +\huno(K_0^\delta)\le
\huno(K_0)$, which is an obvious consequence of
the minimality of  $u_0^\delta$ and of the fact that
$g_0^\delta=g(0)=0$.
\end{proof}

\begin{lemma}\label{Helly2} Let $\lambda$ be the constant in
Lemma~\ref{estim}. Then $K_{\delta}(t)\in\Kml$ and
$K(t)\in \Kml$ for every $t\in[0,T]$.
\end{lemma}
\begin{proof}
By Lemma~\ref{estim} we have $\huno(K_{\delta}(t))\le\lambda$ for every
$t\in[0,T]$ and every $\delta>0$. By Theorem~\ref{Golab2} this
implies $K(t)\in\Kml$ for every $t\in[0,T]$.
\end{proof}

For every $t\in[0,T]$ let $u(t)$ be  a solution of  the minimum problem
(\ref{e}) which defines $\E(g(t),K(t))$.
\begin{lemma}\label{l2}
For every $t\in[0,T]$ we have
$\nabla u_{\delta}(t)\to\nabla u(t)$ strongly in $L^2(\Om;\R^2)$.
\end{lemma}
\begin{proof} As $u_{\delta}(t)$ is a solution of the minimum  problem
(\ref{e}) which defines $\E(g_{\delta}(t),K_{\delta}(t))$, and
$g_\delta(t)\to g(t)$ strongly in $H^1(\Om)$, the
conclusion follows from
Theorem~\ref{convsol} and Lemma \ref{Helly2}.
\end{proof}

The following lemma shows the minimality of the set
$K(t)$ for the functional
$\E(g(t),\cdot)$ with respect to sets $K$ containing $K(t)$.

\begin{lemma}\label{condb}
For every $t\in[0,T]$ we have
\begin{equation}\label{pt}
\E(g(t),K(t))\leq \E(g(t),K)\quad\forall\,K\in\Kmf\,,\ K\supset
K(t)\,.
\end{equation}
\end{lemma}
\begin{proof} Let us fix $t\in[0,T]$ and $K\in\Kmf$ with $K\supset K(t)$.
Since  $K_\delta(t)$ converges
to $K(t)$ in the Hausdorff metric as $\delta\to0$, it is possible to
construct a a sequence
$(K_\delta)$ in $\Kmf$,
converging to $K$ in the Hausdorff metric, such that
$K_\delta\supset K_\delta(t)$ and
$\huno(K_\delta\setmeno K_\delta(t))\to \huno(K\setmeno K(t))$
as $\delta\to0$.
By Lemma \ref{estim} this implies that $\huno(K_\delta)$ is bounded as
$\delta\to0$.
The main difficulty in the construction of $K_\delta$
is the constraint on the number of connected components.
The proof of the details is quite long, but
elementary, and is given in \cite[Lemma~3.5]{DM-Toa}.

Let $v_\delta$ and $v$ be solutions of the minimum problems
(\ref{e}) which define $\E(g_\delta(t),K_\delta)$ and
$\E(g(t),K)$, respectively. By Theorem~\ref{convsol}
$\nabla v_\delta\to\nabla v$ strongly
in $L^2(\Om;\R^2)$.
The minimality of $K_\delta(t)$ expressed by (\ref{pidelta})
gives
$\E(g_\delta(t),K_\delta(t))\leq \E(g_\delta(t),K_\delta)$,
which implies
$\|\nabla u_\delta(t)\|^2\leq\|\nabla v_\delta\|^2+\huno(K_\delta\setmeno
K_\delta(t))$.
Passing to the limit as $\delta\to0$ and using Lemma~\ref{l2}
we get $\|\nabla u(t)\|^2\leq\|\nabla v\|^2+\huno(K\setmeno K(t))$.
Adding $\huno(K(t))$ to both sides we obtain (\ref{pt}).
\end{proof}

We can now pass to the limit in Lemma~\ref{discr}.

\begin{lemma}\label{ineq}
For every $s,\, t$ with $0\leq s<t\leq T$
\begin{equation}\label{diseg}
\|\nabla u(t)\|^2+\huno(K(t))\leq\|\nabla u(s)\|^2+\huno(K(s))+
2\int_s^t(\nabla u(\tau)|\nabla \gdot(\tau))d\tau\,.
\end{equation}
\end{lemma}
\begin{proof}
Let us fix $s,t$ with $0\leq s<t\leq T$. Given $\delta>0$ let $i$ and $j$
be the
integers such that $t_i^\delta\leq s<t_{i+1}^\delta$ and
$t_j^\delta\leq t<t_{j+1}^\delta$.
Let us define $s_\delta:=t_i^\delta$ and $t_\delta:=t_j^\delta$.
Applying Lemma~\ref{discr} we obtain
\begin{equation}\label{ediscr}
\|\nabla u_\delta(t)\|^2+\huno(K_\delta(t)\setmeno K_\delta(s))\leq
\|\nabla u_\delta(s)\|^2+
2\int_{s_{\delta}}^{t_{\delta}}\!\!\!(\nabla u_\delta(\tau)|\nabla
\gdot(\tau))\,d\tau+\omega(\delta)\,,
\end{equation}
with $\omega(\delta)$ converging to zero as $\delta\to0$.
By Lemma~\ref{l2}
for every $\tau\in[0,T]$ we have $\nabla u_{\delta}(\tau)\to \nabla u(\tau)$
strongly in
$L^2(\Om,\R^2)$ as $\delta\to0$, and by Lemma~\ref{estim} we have
$\|\nabla u_\delta(\tau)\|\le \lambda$ for every $\tau\in[0,T]$.

Given $\e>0$, let
$K^{\e}(s):=\{x\in\overline\Om: {\rm dist}(x,K(s))\le \e\}$. As
$K_\delta(s)\subset  K^{\e}(s)$ for $\delta$ small enough, we have
$K_\delta(t)\setmeno K^{\e}(s)\subset  K_\delta(t)\setmeno K_\delta(s)$.
Applying Theorem~\ref{Golab2} with $U=\R^2\setmeno K^{\e}(s)$ we get
$$
\huno(K(t)\setmeno K^{\e}(s)) \le
\liminf_{\delta\to0} \, \huno(K_\delta(t)\setmeno K^{\e}(s))
\le \liminf_{\delta\to0} \, \huno(K_\delta(t)\setmeno K_\delta(s))\,.
$$
Passing to the limit as $\e\to0$ we obtain
$$
\huno(K(t)\setmeno K(s))\leq\liminf_{\delta\to0}\,
\huno(K_\delta(t)\setmeno K_\delta(s))\,.
$$
Passing now to the limit in (\ref{ediscr}) as $\delta\to0$ we
obtain (\ref{diseg}).
\end{proof}

We are now in a position to prove the absolute continuity of the
function $t\mapsto \E(g(t),K(t))$ and to compute its derivative.
\begin{lemma}\label{abscont}
The function $t\mapsto \E(g(t),K(t))$ is absolutely continuous on
$[0,T]$ and
\begin{equation}\label{gtkt}
\frac{d}{dt}\E(g(t),K(t))=2(\nabla u(t)|\nabla\gdot(t))\qquad\hbox{ for
a.e.\ }t\in[0,T]\,.
\end{equation}
\end{lemma}

\begin{proof}Let $0\leq s<t\leq T$.
{}From the previous lemma we get
\begin{equation}\label{alto}
\E(g(t),K(t))-\E(g(s),K(s))\leq2\int_s^t(\nabla
u(\tau)|\nabla\gdot(\tau))\,d\tau\,.
\end{equation}
On the other hand, by Lemma~\ref{condb} we have
$\E(g(s),K(s))\leq \E(g(s),K(t))$. It is easy to see that the Frechet
differential $d\E(g,K)$ of $\E(g,K)$ (with respect to $g$) is given by
\begin{equation}\label{dE}
d\E(g,K)\,h=2\int_{\Om\setminus K}\nabla u_g\nabla h\,dx\,,
\end{equation}
where $u_g$ is a solution of the minimum  problem (\ref{e}) which
defines $\E(g,K)$. Therefore we have
$$
\E(g(t),K(t))-\E(g(s),K(t))=
2\int_s^t(\nabla u(\tau, t)|\nabla \gdot(\tau))\,d\tau\,,
$$
where $u(\tau,t)$ is a solution of the minimum problem~(\ref{e}) which defines
$\E(g(\tau),K(t))$.
Together with the inequality $\E(g(s),K(s))\leq \E(g(s),K(t))$, this implies
\begin{equation}\label{basso}
\E(g(t),K(t))-\E(g(s),K(s))\geq
2\int_s^t(\nabla u(\tau, t)|\nabla \gdot(\tau))\,d\tau\,.
\end{equation}
Since there exists a constant $C$ such that
$\|\nabla u(\tau)\|\leq\|\nabla g(\tau)\|\leq C$ and
$\|\nabla u(\tau,t)\|\leq\|\nabla g(\tau)\|\leq C$ for $s\leq\tau\leq t$,
from (\ref{alto}) and (\ref{basso})
we obtain
$$
\big|\E(g(t),K(t))-\E(g(s),K(s))\big|\leq 2\,C
\int_s^t\|\nabla \gdot(\tau)\|\,d\tau\,,
$$
which proves that the  function $t\mapsto \E(g(t),K(t))$ is absolutely
continuous.

As $\nabla u(\tau,t)\to\nabla u(t)$ strongly in $L^2(\Om;\R^2)$ when
$\tau\to t$, if we divide (\ref{alto}) and (\ref{basso}) by $t-s$, and
take the limit as $s\to t-$ we obtain (\ref{gtkt}).
\end{proof}

The result of the previous lemma can be expressed equivalently in the
following way.

\begin{lemma}\label{der0}
The function $t\mapsto \E(g(t),K(t))$ is absolutely continuous on
$[0,T]$ and
\begin{equation}\label{gtks}
\frac{d}{ds}\E(g(t),K(s))\Big|\lower1.5ex\hbox{$\scriptstyle
s=t$}=0\quad\hbox{for a.e.\ }t\in[0,T]\,.
\end{equation}
\end{lemma}

\begin{proof}
Let $\Theta$ be the set defined in Proposition~\ref{diff3}. By
(\ref{dE}) and by Theorem~\ref{convsol} the differential
$d\E(g,K(s))$ tends to $d\E(g(t),K(t))$ as $(g,s)$ tends to
$(g(t),t)$ in $H^1(\Om)\times\R$. Let us fix a point $t$ in $\Theta$
such that the function $s\mapsto \E(g(s),K(s))$ is differentiable at
$s=t$ and $t$ is a Lebesgue point of $\dot g$.
For every $s\in[0,T]$ we have
\begin{eqnarray*}
&\E(g(s),K(s))- \E(g(t),K(t))=\\
&=\E(g(s),K(s))-\E(g(t),K(s)) + \E(g(t),K(s))-\E(g(t),K(t))=\\
&= 
\displaystyle
\int_t^s d\E(g(\tau),K(s))\,\dot g(\tau)\,d\tau +
\E(g(t),K(s))-\E(g(t),K(t))\,.
\end{eqnarray*}
Dividing by $s-t$ and taking the limit as $s\to t$ we obtain
$$
\frac{d}{ds}\E(g(s),K(s))\Big|\lower1.5ex\hbox{$\scriptstyle s=t$}
=
d\E(g(t),K(t))\,\dot g(t) +
\frac{d}{ds}\E(g(t),K(s))\Big|\lower1.5ex\hbox{$\scriptstyle s=t$}
\,.
$$
The conclusion follows from (\ref{dE}) and Lemma~\ref{abscont}.
\end{proof}

The properties of the function $K\colon[0,T]\to\K$ are summarized by
the following theorem, which is an immediate consequence of
Lemmas~\ref{Helly2}, \ref{condb}, \ref{abscont},
and~\ref{der0}.
\begin{theorem}\label{kt}
Let $m\ge1$, let $g\in AC([0,T];H^1(\Om))$, and let $K_{0}\in\Kmf$. Then
the function $K\colon[0,T]\to\K$ introduced at the end of the last
lecture satisfies the following properties:
\smallskip
\begin{itemize}
\item[(a)] \hfil $\displaystyle \vphantom{\frac{d}{ds}}
K(0)=K_0$, \hfil
\item[(b)] \hfil $\displaystyle \vphantom{\frac{d}{ds}}
K_0\subset K(s)\subset K(t) \,\hbox{ and }\, K(t)\in\Kmf \,
\hbox{ for }\, 0\le s\le t\le T$, \hfil
\item[(c)] \hfil $\displaystyle \vphantom{\frac{d}{ds}}
\hbox{for  }\, 0\le t \le T\quad\E(g(t),K(t))\leq \E(g(t),K)
\quad\forall \, K\in\Kmf,\,\  K\supset K(t)$,\hfil
\item[(d)]\hfil $\displaystyle \vphantom{\frac{d}{ds}}
t\mapsto \E(g(t),K(t)) \hbox{ is
absolutely continuous on }[0,T]$, \hfil
\item[(e)]\hfil$\displaystyle\frac{d}{ds}\E(g(t),K(s))
\Big|\lower1.5ex\hbox{$\scriptstyle s=t$}=0\quad
\hbox{for a.e.\ }t\in[0,T]$.\hfil
\end{itemize}
\smallskip
Moreover every function $K\colon[0,T]\to\Kmf$ which satisfies (a)--(e)
satisfies also
\begin{itemize}
\item[(f)] \hfil$\displaystyle\frac{d}{dt}\E(g(t),K(t))=
2(\nabla u(t)|\nabla \gdot(t))\quad \hbox{ for a.e.\ } t\in[0,T]$,\hfil
\end{itemize}
where $u(t)$ is a solution of the minimum  problem (\ref{e}) which defines
$\E(g(t),K(t))$.
\end{theorem}

In our continuous-time model, the function $K\colon[0,T]\to\Kmf$
represents the quasi-static irreversible evolution of the crack
starting from $K_0$ (condition (a)) under
the action of the boundary displacement $g(t)$. Condition (b)
reflects the {\it irreversibility of the evolution\/}
and the {\it absence of a healing process\/}. Condition (c) is a
unilateral minimality condition.
Condition (e) says that, for
almost every $t\in[0,T]$,
the total energy $s\mapsto \E(g(t) , K(s))$ is
stationary at ${s=t}$. Conditions (c) and (e) together lead to 
Griffith's analysis
of the energy balance in our model, and, under some very mild
regularity assumptions on the cracks, allow to express the classical
Griffith's criterion for crack growth in terms of the stress intensity
factors at the tips of the cracks (see \cite[Section~8]{DM-Toa}).

We underline that, although we can not exclude that the surface energy
$\huno(K(t))$ may present some jump discontinuities in time
(see \cite[Section 4.3]{FraMar3}), in our result
{\it the total energy is always an
absolutely continuous function of time\/} by condition~(d).

If $\partial_D\Om$ is sufficiently smooth, we can integrate  by parts
the right hand side of (f) and, taking into account the Euler equation
satisfied by $u(t)$, we obtain
\begin{equation}\label{energy}
\frac{d}{dt}\E(g(t),K(t))=
2\int_{\partial_D\Om\setminus K(t)}
\frac{\partial u(t)}{\partial \nu}\, \gdot(t)\,d\huno
\quad \hbox{ for a.e.\ } t\in[0,T]\,,
\end{equation}
where $\nu$ is the outer unit normal to $\partial\Om$. Since the
right hand side of (\ref{energy}) is the power of the force exerted
on the boundary to obtain the displacement $g(t)$ on
$\partial_D\Om\setmeno K(t)$,
equality (\ref{energy}) expresses the {\it conservation of energy\/}
in our quasi-static model, where all kinetic effects are neglected.

The discrete-time model described in the previous lecture
turns out to be a useful tool for the proof of the
existence of a solution $K(t)$ of the problem considered in
Theorem~\ref{kt}, and provides also an effective way for the numerical
approximation of this solution (see \cite{BFM}),
since many algorithms have been
developed for the numerical solution of minimum problems of the form
(\ref{pidelta}) (see, e.g., \cite{BZ}, \cite{Rich}, \cite{RM}, \cite{Bou},
\cite{Ch}, \cite{Bou-Cha}).
\end{section}


\end{document}